\documentclass{amsart}
\usepackage{graphics}
\newtheorem{thm}{Theorem}
\newtheorem{lem}[thm]{Lemma}
\newtheorem{sublem}{Sublemma}

\newtheorem{conj}{Conjecture}
\begin{document}
\title[3-superbridge knots]{There are only finitely many 3-superbridge knots}
\author[C.B.\ Jeon and G.T.\ Jin]{Choon Bae Jeon and Gyo Taek Jin}
\address{Korea Advanced Institute of Science and Technology}
\email{jcbout@math.kaist.ac.kr, trefoil@kaist.ac.kr}
\begin{abstract}
Although there are infinitely many knots with superbridge index
$n$ for every even integer $n\ge4$, there are only finitely many
knots with superbridge index 3.
\end{abstract}
\maketitle

\vskip0.3in
\section{Introduction}
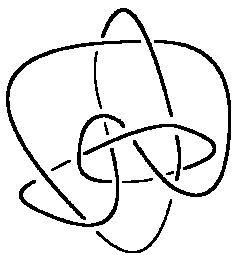
\begin{figure}[b]
\begin{picture}(75, 90)(3, 0)
\thicklines \qbezier(34,80)(45,110)(60,50)
\qbezier(62,42)(64,34)(62,26) \thinlines
\qbezier(60,18)(52,-14)(32,6) \thicklines
\qbezier(27,11)(-40,78)(48,78) \qbezier(54,78)(82,78)(82,50)
\qbezier(82,50)(82,12)(62,22) \qbezier(62,22)(56,25)(46,40)
\qbezier(42,46)(40,49)(35,50) \qbezier(35,50)(25,50)(25,33)
\qbezier(25,33)(25,25)(37,25) \thinlines
\qbezier(42,25)(50,25)(53,26)
\qbezier(58.7,27.9)(59,28)(60.5,28.5)
\thicklines\qbezier(65,30)(89,38)(62,46)
\qbezier(62,46)(52,48)(27,35.5) \thinlines
\qbezier(21,32.5)(24,34)(15,29.5) \thicklines
\qbezier(8,25)(-2,20)(8,15) \qbezier(8,15)(48,-5)(37.5,37)
\thinlines \qbezier(36,43)(36,43)(35,47)
\qbezier(33.5,53)(30,65)(32,74)
\end{picture}
\caption{Bridges---maximal overpasses}\label{fig:overpass}
\end{figure}

Throughout this article a {\em knot\/} is a piecewise smooth
simple closed curve embedded in the three dimensional Euclidean
space $\mathbb R^3$. For a knot $K$, its equivalence class, under
piecewise smooth homeomorphisms of $\mathbb R^3$ mapping one knot
onto another, will be referred to as the {\em knot type\/} of $K$
and denoted by $[K]$.

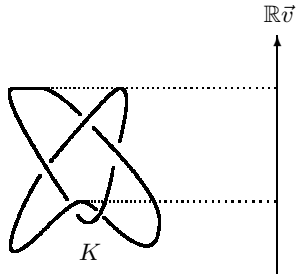
\begin{figure}[t]
\small
\begin{picture}(105,100)(20,45)
\put(120,40){\vector(0,1){90}} \put(115,135){$\mathbb R\vec v$}
{\thicklines \qbezier(40,68)(12,110)(20,110)
\qbezier(20,110)(25,110)(30,110) \qbezier(30,110)(35,110)(45,100)
\qbezier(50,95)(75,70)(75,60) \qbezier(75,60)(75,40)(54,61)
\qbezier(51,64)(48,67)(45,67) \qbezier(45,67)(42,67)(32,57)
\qbezier(32,57)(12,37)(22,62) \qbezier(22,62)(26,71)(30,77)
\qbezier(34,82.4)(57,110)(60,110) \qbezier(60,110)(65,110)(60,90)
\qbezier(57.5,80)(54,59)(48,59) \qbezier(48,59)(46,59)(44,62) }
\qbezier[45](30,110)(65,110)(119,110)
\qbezier[38](45,67)(65,67)(119,67) \put(45,45){$K$}
\end{picture}
\caption{$b_{\vec v}(K)=3$}\label{fig:local-max}
\end{figure}

In 1954, Schubert introduced the {\em bridge index\/} of
knots~\cite{schubert}. In a knot diagram, maximal overpasses are
called {\em bridges}. Figure~\ref{fig:overpass} shows a knot
diagram with seven bridges which are drawn with thick arcs. The
bridge index of a knot is defined to be the minimum number of
bridges in all the possible diagrams of knots in its knot type.
An equivalent definition can be given in the following way. Given
a knot $K$ and a unit vector
 $\vec v$ in $\mathbb R^3$, we define $b_{\vec v}(K)$ as the number of
connected components of the preimage of the set of local maximum
values of the orthogonal projection $K\to\mathbb R\vec v$.
Figure~\ref{fig:local-max} illustrates an example. The {\em
bridge number\/} of $K$ is defined by the formula
$$b(K)=\min_{\|\vec v\|=1}b_{\vec v}(K).$$
It is known that the bridge index can be defined by the formula
$$b[K]=\min_{K^\prime\in[K]}b(K^\prime)
=\min_{K^\prime\in[K]}\min_{\|\vec v\|=1}b_{\vec v}(K^\prime).$$

In 1987, Kuiper modified the alternative definition of bridge
index to define another knot invariant called {\em superbridge
index\/}~\cite{kuiper}. Given a knot $K$, the {\em superbridge
number\/} of $K$ is defined by $$s(K)=\max_{\|\vec v\|=1}b_{\vec
v}(K)$$ and the {\em superbridge index\/} of $K$ by
$$s[K]=\min_{K^\prime\in[K]}s(K^\prime)=\min_{K^\prime\in[K]}\max_{\|\vec
v\|=1} b_{\vec v}(K^\prime).$$ He used Milnor's total
curvature~\cite{milnor} to prove that any nontrivial knot $K$
satisfies the inequality:
\begin{equation}
b[K]<s[K]\label{eq:b<sb}
\end{equation}
He computed the superbridge index for all torus knots.
\begin{thm}[Kuiper]\label{thm:kuiper-torus}
For any two coprime integers $p$ and $q$, satisfying $2\le p<q$,
the superbridge index of the torus knot of type $(p,q)$ is
$\min\{2p, q\}$.
\end{thm}

\section{Odd-superbridge knots}
As the knots having bridge index $n$ are referred to as
$n$-bridge knots, we will call the knots with superbridge index
$n$ as {\em $n$-superbridge knots}.

Because nontrivial knots have bridge index at least 2, the
inequality~(\ref{eq:b<sb}) implies that nontrivial knots have
superbridge index at least $3$. By the same reason, 3-superbridge
knots are 2-bridge knots, in particular, prime knots. According
to Theorem~\ref{thm:kuiper-torus}, trefoil knot is the only torus
knot with superbridge index 3. Figure eight knot is also a
3-superbridge knot~\cite{jin-poly,trautwein}. No other
3-superbridge knots are known yet. Our main theorem asserts that
there are only finitely many 3-superbridge knots.

\begin{thm}\label{thm:3-sup} There are only finitely many 3-superbridge knots.
\end{thm}

By Theorem~\ref{thm:kuiper-torus}, we know that the torus knot of
type $(n,nk+1)$ has superbridge index $2n$, for $n\ge2$ and
$k\ge2$. Therefore, for any even number $2n\ge4$, there are
infinitely many $2n$-superbridge knots.

Because it is natural to expect that more knotting would increase
the superbridge number, we expect that there are infinitely many
$n$-superbridge knots for any positive integer $n\ge4$.

\begin{conj}\label{conj:infty-odd-supbr}
There are infinitely many $(2n-1)$-superbridge knots for any
positive integer $n\ge3$.
\end{conj}

A result in \cite{jin-supbr} implies that a connected sum of any
torus knot $K$ with a trefoil knot $T$ satisfies the inequality
$s[K\sharp T]\le s[K]+1$. As it is generally expected that the
superbridge index of a composite knot would be bigger than that
of any of the factor knots, which is true for bridge index, we
make Conjecture~\ref{conj:add-tref} which implies
Conjecture~\ref{conj:infty-odd-supbr}. This conjecture is valid if
$K$ is a trefoil knot or figure eight knot.

\begin{conj}\label{conj:add-tref}
Every nontrivial knot $K$ satisfies $$s[K\sharp T]=s[K]+1$$ where
$T$ is a trefoil knot.
\end{conj}

By Theorem~\ref{thm:kuiper-torus}, we know that the torus knot of
type $(n,2n-1)$ is a $(2n-1)$-superbridge knot. This knot is the
closure of the $n$-braid
$(\sigma_1\sigma_2\cdots\sigma_{n-1})^{2n-1}$. On the other hand,
the torus knot of type $(2,2k-1)$ is the closure of the 2-braid
$\sigma_1^{2k-1}$ and has superbridge index $4$ if $k\ge3$. For
these torus knots, inserting a full twist $\sigma_1^2$ does not
increase the superbridge index. This fact encourages us to
consider Conjecture~\ref{conj:add-twists} which also implies
Conjecture~\ref{conj:infty-odd-supbr}.

\begin{conj}\label{conj:add-twists}
For $n\ge 3$ and $k\ge0$, the closure of the  $n$-braid
\begin{equation}\label{eqn:twisted-torus}
\sigma_{1}^{2k}(\sigma_1\sigma_2\cdots\sigma_{n-1})^{2n-1}
\end{equation}
is a $(2n-1)$-superbridge knot.
\end{conj}

A theorem of Stallings~\cite{stallings} implies that the closure
of the braid in (\ref{eqn:twisted-torus}) for $k\ge0$ is a fibred
knot with the fibre surface obtained by Seifert's algorithm on the
closed braid diagram. This surface is the one with minimal genus,
which is $(n-1)^2+k$. Therefore for each $n$, such knots are all
distinct. Notice that the braid (\ref{eqn:twisted-torus}) is
positive and the diagram of its closure is visually prime.
According to Cromwell~\cite{cromwell}, they are all prime knots.
The primeness of these knots makes
Conjecture~\ref{conj:add-twists} more interesting than
Conjecture~\ref{conj:add-tref}.

The second author would like to thank Paul Melvin for a
discussion which inspired Conjecture~\ref{conj:add-tref} and to
thank Dale Rolfsen for bringing Stallings' theorem to his
attention.

\section{Proof of Theorem~\ref{thm:3-sup}}
Our proof of Theorem~\ref{thm:3-sup} requires two main tools. The
first is Lemma~\ref{lem:straighten} and the second is {\em
quadrisecant\/} which is a straight line intersecting a knot at
four distinct points. According to~\cite{morton-mond,pannwitz},
every nontrivial knot has a quadrisecant.

\begin{lem}\label{lem:straighten}
Given a knot $K$, let $K^\prime$ be a knot obtained by replacing
a subarc of $K$ with a straight line segment joining the end
points of the subarc. Then $s(K)\ge s(K^\prime)$.
\end{lem}

{\sc Proof:} Given a unit vector $\vec v$, let
$g\colon(-1,2)\to\mathbb R\vec v$ be a parametrization of the
orthogonal projection of an open neighborhood of the subarc into
$\mathbb R\vec v$, where the subarc corresponds to the closed
interval $[0,1]$. Then the projection of a neighborhood of the
straight line segment in $K^\prime$ can be parametrized by
$$g^\prime(t)=\begin{cases}(1-t)g(0)+t g(1)&\mbox{if } t\in[0,1]\\
                          g(t)&\mbox{if }t\in(-1,0]\cup[1,2).
              \end{cases}
$$ Since $g^\prime$ has no more local maxima than $g$,  we have
$b_{\vec v}(K)\ge b_{\vec v}(K^\prime)$ for any $\vec v$.
Therefore $s(K)\ge s(K^\prime)$.~\qed

\medskip
Let $K$ be a 3-superbridge knot with superbridge number 3, namely,
$s[K]=s(K)=3$, and let $\mathcal Q$ be a quadrisecant of $K$.
Then $K-\mathcal Q$ consists of four disjoint open arcs
$l_1,l_2,l_3$ and $l_4$.  Let $\bar l_i$ and $\tilde l_i$ denote
$\pi(l_i)$ and $\pi(\mathcal Q\cup l_i)$, respectively, where
$\pi\colon\mathbb R^3\to \mathcal Q^\perp$ is the orthogonal
projection of $\mathbb R^3$ onto a plane $\mathcal Q^\perp$
perpendicular to the quadrisecant. Applying
Lemma~\ref{lem:straighten} wherever needed, we may assume that
the only singular points of $\pi(K)$ are a set of finitely many
transversal double points together with a quadruple point
$\pi(\mathcal Q)$. For every open subarc $l$ of $K$, write
$b_{\vec v}(K\mid l)$ for the number of local maxima of
$K\to\mathbb R\vec v$ on $l$. Since each $\tilde l_i$ is a closed
loop in $\mathcal Q^\perp$, we must have
\begin{equation}\label{eqn:any-v}
b_{\vec v}(K\mid l_i)\ge1 \mbox{\quad or\quad} b_{-\vec v}(K\mid
l_i)\ge1
\end{equation}
for every unit vector $\vec v\in \mathcal Q^\perp$.

For a straight line $\rho$ in $\mathcal Q^\perp$, let $\vec
v_\rho$ denote a unit vector in $\mathcal Q^\perp$ perpendicular
to $\rho$.

\begin{sublem}\label{sublem:no-self-xing}
We may assume that $\bar l_i$ has no self-crossings, for each
$i=1,2,3,4$.
\end{sublem}

{\sc Proof:} Suppose $\bar l_i$ has a self-crossing. The we can
choose a loop $\lambda$ of $\bar l_i$ which is {\em minimal\/} in
the sense that no proper subarc of $\lambda$ is another loop.
Then $\lambda$ bounds an open disk $\delta$ in $\mathcal Q^\perp$.

If $\pi(K)\cap\delta=\emptyset$, we can eliminate this loop
together with its crossing by a move as described in
Lemma~\ref{lem:straighten} without changing the knot type.

If $\bar l_i$ passes through $\delta$, then among the half-lines
starting from $\pi(\mathcal Q)$ and passing through $\delta$, we
are able to find one, say $\rho$, which meets $\bar l_i$ at least
three times. Then we have
\begin{equation}\label{eqn:2-2}
b_{ \vec v_\rho}(K\mid l_i)\ge2\quad\text{and}\quad b_{-\vec
v_\rho}(K\mid l_i)\ge2.
\end{equation}
This, together with the fact (\ref{eqn:any-v}), implies
\begin{equation}\label{eqn:4-4}
b_{\vec v_\rho}(K)\ge\sum_{1\le j\le 4}b_{\vec v_\rho}(K\mid
l_j)\ge 4 \quad \text{or}\quad b_{-\vec v_\rho}(K)\ge\sum_{1\le
j\le 4}b_{-\vec v_\rho}(K\mid l_j)\ge 4
\end{equation}
which contradicts $s(K)=3$.

If $\bar l_j$ passes through $\delta$, for some $j\ne i$, then
among the half-lines starting from $\pi(\mathcal Q)$ and passing
through $\delta$, we are able to find one, say $\rho$, which
crosses $\bar l_j$. Then, for $\vec w=\vec v_\rho$ or $\vec
w=-\vec v_\rho$, we have
\begin{equation}\label{eqn:2111}
b_{\vec w}(K\mid l_i)\ge2,\quad b_{-\vec w}(K\mid l_i)\ge1,\quad
b_{\vec w}(K\mid l_j)\ge1,\quad b_{-\vec w}(K\mid l_j)\ge1.
\end{equation}
This, together with the fact (\ref{eqn:any-v}), implies
(\ref{eqn:4-4}) which contradicts $s(K)=3$.~\qed

\begin{sublem}\label{sublem:star-shaped}
We may assume that each $\tilde l_i$ bounds an open disk
$\delta_i$ in $\mathcal Q^\perp$ which is star-shaped with
respect to $\pi(\mathcal Q)$.
\end{sublem}

{\sc Proof:} By Sublemma~\ref{sublem:no-self-xing}, we know that
$\tilde l_i$ bounds an open disk $\delta_i$ in $\mathcal Q^\perp$.
If $\delta_i$ is not star-shaped, there exists a half-line $\rho$
in $\mathcal Q^\perp$ starting at $\pi(\mathcal Q)$ and meeting
$\bar l_i$ more than once. If $\rho$ meets $\bar l_i$ at three or
more points, then the condition (\ref{eqn:2-2}) holds. Therefore
we reach the same contradiction as in (\ref{eqn:4-4}). Suppose
$\rho$ meets $\bar l_i$ at two points. Then there exist two open
disks $R$ and $S$ bounded by $\rho$ and $\bar l_i$ as in
Figure~\ref{fig:meet-at-2}(i). If $\bar l_j$ meets $R\cup S$,
there is a half line $\rho^\prime$ starting from $\pi(\mathcal
Q)$ crossing $\bar l_j$ at a point in $R\cup S$ as indicated by
Figure~\ref{fig:meet-at-2}(ii). Then, for $\vec w=\vec
v_{\rho^\prime}$ or $\vec w=-\vec v_{\rho^\prime}$, the condition
(\ref{eqn:2111}) holds. This leads to the same contradiction as
in (\ref{eqn:4-4}). If there are no arcs of $\pi(K)$ inside
$R\cup S$, we can straighten a part of $l_i$ as in
Figure~\ref{fig:meet-at-2}(iii) without changing the knot type.
By Lemma~\ref{lem:straighten}, this move doesn't increase the
superbridge number. Since $s[K]=s(K)=3$, this move cannot
decrease the superbridge number either. Only finitely many of
such modifications are necessary to deform $\delta_i$ into a
region star-shaped with respect to $\pi(\mathcal Q)$.~\qed

\begin{figure}[t]
\small
\begin{center}
\begin{picture}(85,44)(-5,-19)
{\thicklines \qbezier(0,0)(28,-24)(40,0)
\qbezier(40,0)(70,60)(60,-5) \qbezier(60,-5)(52,-57)(0,0) }
\put(0,0){\line(1,0){80}} \put(72,-8){$\rho$}
\put(20,-9){$R$}\put(51,4){$S$} \put(42,-23){$\bar l_i$}
\put(-5,-28){(i)}
\end{picture}
\qquad
\begin{picture}(85,14)(-5,-19)
{\thicklines \qbezier(0,0)(28,-24)(40,0)
\qbezier(40,0)(70,60)(60,-5) \qbezier(60,-5)(52,-57)(0,0)
\qbezier(44,25)(64,3)(74,3) } \put(0,0){\line(1,0){80}}
\qbezier(0,0)(0,0)(76,19) \put(72,-8){$\rho$}
\put(73,12){$\rho^\prime$} \put(0,0){\vector(-1,4){3.5}}
\put(1,9){$\vec v_{\rho^\prime}$} \put(42,-23){$\bar l_i$}
\put(36,20){$\bar l_j$} \put(-5,-28){(ii)}
\end{picture}
\qquad
\begin{picture}(85,44)(-5,-19)
\qbezier[25](0,0)(28,-24)(40,0) \qbezier[55](40,0)(70,60)(60,-5)
{\thicklines \qbezier(0,0)(0,0)(60,-5)
\qbezier(60,-5)(52,-57)(0,0) } \put(0,0){\line(1,0){80}}
\put(72,-8){$\rho$} \put(-5,-28){(iii)}
\end{picture}
\end{center}
\caption{}\label{fig:meet-at-2}
\end{figure}

\begin{sublem}\label{sublem:delta-delta}
None of the following conditions hold when $h,i,j,k$ are distinct
elements of $\{1,2,3,4\}$.
\begin{eqnarray}
& \delta_i\cap\delta_j\cap\delta_k\ne\emptyset
      \label{eqn:ijk}\\
& \delta_i\cap\delta_j\ne\emptyset\ \text{and}\
\delta_h\cap\delta_k\ne\emptyset
      \label{eqn:ij-hk}\\
& \delta_i\cap\delta_j\ne\emptyset,\
\delta_i\cap\delta_k\ne\emptyset\ \text{and}\
\delta_i\cap\delta_h\ne\emptyset
      \label{eqn:ij-ik-ih}
\end{eqnarray}
\end{sublem}

{\sc Proof:} For each of the three conditions, we will choose a
line $\rho$ in $\mathcal Q^\perp$ which meets $\pi(K)$ at least
eight times. Then we must have $b_{\vec v_\rho}(K)\ge4$, which
contradicts $s(K)=3$.

\noindent{\sc Condition (\ref{eqn:ijk}):} If this condition is
true, we can choose two points
$P_1\in\delta_i\cap\delta_j\cap\delta_k$ and $P_2\in\delta_h$ so
that the straight line $\rho$ joining $P_1$ and $P_2$ does not
pass through $\pi(\mathcal Q)$. Then each $\bar l_a$ crosses
$\rho$ at least twice, for $a=h,i,j,k$.

\noindent{\sc Condition (\ref{eqn:ij-hk}):} If this condition is
true, we can choose two points $P_1\in\delta_i\cap\delta_j$ and
$P_2\in\delta_h\cap\delta_k$ so that the straight line $\rho$
joining $P_1$ and $P_2$ does not pass through $\pi(\mathcal Q)$.
Then again, each $\bar l_a$ crosses $\rho$ at least twice, for
$a=h,i,j,k$.

\noindent{\sc Condition (\ref{eqn:ij-ik-ih}):} If this condition
is true, we can choose three points $P_a\in\delta_i\cap\delta_a$
for $a=j,k,h$, so that the three straight lines determined by
pairs of $P_a$'s do not pass through $\pi(\mathcal Q)$. Since
(\ref{eqn:ijk}) cannot occur, every edge of the triangle
$\triangle P_jP_kP_h$ meets $\pi(K)$ in even number of times.
There are two subcases to consider:

\noindent{\sc Subcase (\ref{eqn:ij-ik-ih}.1):} If $\pi(\mathcal
Q)$ is contained inside $\triangle P_jP_kP_h$, the boundary of
this triangle meets  $\pi(K)$ at least eight times. Therefore
there is an edge, say $\overline{P_jP_k}$, meeting  $\pi(K)$ at
least four times. Let $\rho$ be the extension of $P_jP_k$.

\noindent{\sc Subcase (\ref{eqn:ij-ik-ih}.2):} If $\pi(\mathcal
Q)$ is contained outside of $\triangle P_jP_kP_h$, there is one
vertex of  $\triangle P_jP_kP_h$, say $P_h$, such that the
straight line segment joining $P_h$ and $\pi(\mathcal Q)$ crosses
the edge $\overline{P_jP_k}$. Since (\ref{eqn:ijk}) cannot occur,
the edge $\overline{P_jP_k}$ meets  $\pi(K)$ at least four times.
Let $\rho$ be the extension of $P_jP_k$.

For the above two subcases, $\rho$ meets $\pi(K)$ at least twice
on either side of the extension.

Consequently, $\rho$ meets $\pi(K)$ at least eight times as
required.~\qed

\begin{sublem}\label{sublem:half-twists}
We may assume that the only crossings in $\bar l_i\cup\bar l_j$
are a set of finitely many consecutive half twists.
\end{sublem}

{\sc Proof:} It follows from
Sublemmas~\ref{sublem:no-self-xing}--\ref{sublem:delta-delta}
that if a half-line starting from $\pi(\mathcal Q)$ meets $\bar
l_i$ and $\bar l_j$, then it meets each of them exactly once and
no other $\bar l_a$'s. Notice that the region
$\delta_i\cap\delta_j$ is also an open disk which is star-shaped
with respect to $\pi(\mathcal Q)$ such that the only singular
points along $\partial(\delta_i\cap\delta_j)$ are the quadruple
point $\pi(\mathcal Q)$ and the half-twists.~\qed

\begin{figure}[h]
\small
\begin{picture}( 75, 90)(-37.5,-18.75)
\put(-15,67.5){\line(1,0){30}} \put( 15,67.5){\line(0,-1){22.5}}
\put( 15,45){\line(-1,0){30}} \put(-15,45){\line(0,1){22.5}}
\qbezier(15,60)(37.5,60)(37.5,37.5)
\qbezier(37.5,37.5)(37.5,15)(0,0)
\qbezier(0,0)(-18.75,-7.5)(-16.875,-9.375)
\qbezier(-16.875,-9.375)(-15,-11.25)(0,0)
\qbezier(0,0)(30,22.5)(30,37.5)
\qbezier(30,37.5)(30,52.5)(15,52.5)
\qbezier(-15,52.5)(-30,52.5)(-30,37.5)
\qbezier(-30,37.5)(-30,22.5)(0,0)
\qbezier(0,0)(15,-11.25)(16.875,-9.375)
\qbezier(16.875,-9.375)(18.75,-7.5)(0,0)
\qbezier(0,0)(-37.5,15)(-37.5,37.5)
\qbezier(-37.5,37.5)(-37.5,60)(-15,60) \put(-5,-26){(a)}
\end{picture}
\quad\quad
\begin{picture}(135, 97.5)(-67.5,-18.75)
\put(-6,-26){(b)} \put(-45,67.5){\line(1,0){30}}
\put(-15,67.5){\line(0,-1){22.5}} \put(-15,45){\line(-1,0){30}}
\put(-45,45){\line(0,1){22.5}} \put(15,67.5){\line(1,0){30}}
\put(45,67.5){\line(0,-1){22.5}} \put(45,45){\line(-1,0){30}}
\put(15,45){\line(0,1){22.5}} \put(-15,60){\line(1,0){30}}
\qbezier(45,60)(67.5,60)(67.5,37.5)
\qbezier(67.5,37.5)(67.5,7.5)(0,0)
\qbezier(0,0)(-60,22.5)(-60,37.5)
\qbezier(-60,37.5)(-60,52.5)(-45,52.5)
\qbezier(-15,52.5)(-7.5,52.5)(0,0) \qbezier(0,0)(7.5,-15)(0,-15)
\qbezier(0,-15)(-7.5,-15)(0,0) \qbezier(0,0)(7.5,52.5)(15,52.5)
\qbezier(45,52.5)(60,52.5)(60,37.5)
\qbezier(60,37.5)(60,22.5)(0,0)
\qbezier(0,0)(-67.5,7.5)(-67.5,37.5)
\qbezier(-67.5,37.5)(-67.5,60)(-45,60)
\end{picture}
\caption{}\label{fig:only-2-kinds}
\end{figure}

\medskip
The next sublemma easily follows from
Sublemmas~\ref{sublem:star-shaped}--\ref{sublem:half-twists}.

\begin{sublem}\label{sublem:only-2-kinds}
We may assume that $\pi(K)$ is as in Figure~\ref{fig:only-2-kinds}
up to planar isotopies of $\mathcal Q^\perp$, where each
rectangle contains a pair of parallel arcs or a pair of arcs with
finitely many half-twists.
\end{sublem}

Suppose that $\delta_i\cap\delta_j\ne\emptyset$ and that none of
$\delta_i$ and $\delta_j$ contains the other completely. Consider
a connected component $\delta$ of $\delta_i-(\delta_j\cup l_j)$
which does not meet $\pi(K)$. By
Sublemma~\ref{sublem:half-twists}, we easily see that
$\partial\delta$ has only two singular points of $\pi(K)$. Such a
region will be referred to as a {\em crescent} and the two
singular points the {\em ends} of the crescent. It is possible
for a crescent to have the quadruple point $\pi(\mathcal Q)$ as
one of its ends. In this case, it is possible to have a {\em
loop-crescent\/} which is bounded by one loop which passes
through $\pi(\mathcal Q)$ and is the projection of one subarc of
$K$.

\begin{sublem}\label{sublem:crescent}
We may assume that every crescent which is not a loop-crescent is
alternating, in the sense that, if one of the two arcs on the
boundary of the crescent passes over the other at one end then it
passes under the other at the other end. We may also assume that
no crescent is a loop-crescent.
\end{sublem}
\begin{figure}[h]
\begin{picture}(90,35)(0,-3)
\thicklines \qbezier(0,0)(15,20)(30,10)
\qbezier(30,10)(45,0)(60,10) \qbezier(60,10)(75,20)(90,0)
\qbezier(0,-10)(10,-20)(20,0) \qbezier(20,0)(45,50)(70,0)
\qbezier(70,0)(80,-20)(90,-10) \thinlines
\qbezier[30](30,10)(45,10)(60,10) \qbezier[50](20,0)(45,0)(70,0)
\put(30,10){\circle*{2}} \put(60,10){\circle*{2}}
\put(20,0){\circle*{2}} \put(70,0){\circle*{2}}
\end{picture}
\quad \raise10pt\hbox{$\Longrightarrow$} \quad
\begin{picture}(90,35)(0,-3)
\thicklines \qbezier(0,0)(15,20)(30,10)
\qbezier(60,10)(75,20)(90,0) \put(30,10){\line(1,0){30}}
\qbezier(0,-10)(10,-20)(20,0) \qbezier(70,0)(80,-20)(90,-10)
\put(20,0){\line(1,0){50}} \thinlines
\qbezier[30](30,10)(45,0)(60,10) \qbezier[50](20,0)(45,50)(70,0)
\put(30,10){\circle*{2}} \put(60,10){\circle*{2}}
\put(20,0){\circle*{2}} \put(70,0){\circle*{2}}
\end{picture}
\caption{Non-alternating crescent}\label{fig:reidemeister-2}
\end{figure}
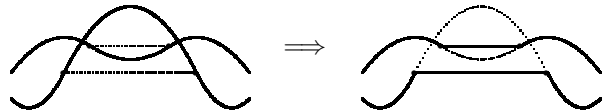

{\sc Proof:}  If a crescent which is away from $\pi(\mathcal Q)$
is non-alternating, we can remove the two crossings at its ends
by straightening two arcs as shown in
Figure~\ref{fig:reidemeister-2}. Since this is a second
Reidemeister move, the knot type does not change. Again by
Lemma~\ref{lem:straighten}, the superbridge number is unchanged.
For a non-alternating crescent whose one end is at $\pi(\mathcal
Q)$, a similar process eliminates the crossing at the other end
if it is not a loop-crescent.

If there is a loop-crescent, we may straighten a small subarc of
the loop near $\pi(\mathcal Q)$ without changing the knot type
and the superbridge number. Then $\mathcal Q$ becomes a {\em
trisecant\/} of the new knot again denoted by $K$. Since $K$ is
nontrivial, the new knot must be obtained from
Figure~\ref{fig:only-2-kinds}(b), and hence its projection must
be as in Figure~\ref{fig:trisec-proj}(a), after straightening out
any unnecessary half-twists. On the other hand, the projection of
a cylindrical neighborhood of a trisecant of an arbitrary knot
has five possible patterns as shown in
Figure~\ref{fig:trisec-proj}(b) where all three arcs are smooth,
up to small perturbations and planar isotopies. All the
combinations of Figure~\ref{fig:trisec-proj}(a) and one of
Figure~\ref{fig:trisec-proj}(b) making a nontrivial knot give
rise to torus knots. Because trefoil knots are the only
3-superbridge torus knots, we can exclude the case of
loop-crescents.~\qed

\begin{figure}[h]
\small
\begin{picture}( 75, 90)(-37.5,-18.75)
\put(-15,67.5){\line(1,0){30}} \put( 15,67.5){\line(0,-1){22.5}}
\put( 15,45){\line(-1,0){30}} \put(-15,45){\line(0,1){22.5}}
\qbezier(15,60)(37.5,60)(37.5,37.5)
\qbezier(37.5,37.5)(37.5,15)(0,0) \qbezier(0,0)(7.5,-15)(0,-15)
\qbezier(0,-15)(-7.5,-15)(0,0) \qbezier(0,0)(30,22.5)(30,37.5)
\qbezier(30,37.5)(30,52.5)(15,52.5)
\qbezier(-15,52.5)(-30,52.5)(-30,37.5)
\qbezier(-30,37.5)(-30,22.5)(0,0)
\qbezier(0,0)(-37.5,15)(-37.5,37.5)
\qbezier(-37.5,37.5)(-37.5,60)(-15,60) \put(-5,-26){(a)}
\end{picture}
\quad
\begin{picture}(160,90)(-50,-6.25)
\put(0,60){\circle{40}} \qbezier(20,60)(0,60)(10,77.32)
\qbezier(-10,77.32)(0,60)(-20,60)
\qbezier(-10,42.68)(0,60)(10,42.68)

\put(50,60){\circle{40}} \qbezier(60,77.32)(50,60)(40,77.32)
\put(70,60){\line(-1,0){40}} \qbezier(40,42.68)(50,60)(60,42.68)

\put(100,60){\circle{40}} \qbezier(110,77.32)(100,60)(90,77.32)
\qbezier(120,60)(100,60)(90,42.68)
\qbezier(80,60)(100,60)(110,42.68)

\put(25,16.70){\circle{40}} \qbezier(35,34.02)(25,16.70)(35,-0.62)
\put(45,16.70){\line(-1,0){40}}
\qbezier(15,34.02)(25,16.70)(15,-0.62)

\put(75,16.70){\circle{40}} \qbezier(85,34.02)(75,16.70)(65,-0.62)
\put(95,16.70){\line(-1,0){40}}
\qbezier(65,34.02)(75,16.70)(85,-0.62) \put(44,-13){(b)}

\end{picture}
\caption{}\label{fig:trisec-proj}
\end{figure}

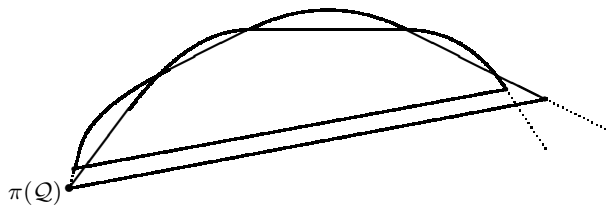
\begin{figure}[b]
\small
\begin{picture}(202.5,75)
\thicklines \put(0,0){\line(3,4){22.5}}
\qbezier(22.5,30)(45,60)(67.5,60) \put(67.5,60){\line(1,0){60}}
\qbezier(127.5,60)(150,60)(165,37.5)
\qbezier(165,37.5)(83.4375,22.5)(1.875,7.5)
\qbezier(1.875,7.5)(1.875,7.5)(3.75,15)
\qbezier(3.75,15)(7.5,30)(37.5,45) \put(37.5,45){\line(2,1){30}}
\qbezier(67.5,60)(97.5,75)(127.5,60)
\put(127.5,60){\line(2,-1){52.5}}
\qbezier(180,33.75)(90,16.875)(0,0) \thinlines
\qbezier[5](0,0)(0.9375,3.75)(1.875,7.5)
\qbezier[14](165,37.5)(172.5,26.25)(180,15)
\qbezier[14](180,33.75)(191.25,28.125)(202.5,22.5)
\put(-23,-4){$\pi(\mathcal Q)$} \put(0,0){\circle*{3}}
\put(1.875,7.5){\circle*{2}} \put(165,37.5){\circle*{2}}
\put(180,33.75){\circle*{2}}
\end{picture}
\caption{$\bar l_i$ and $\bar l_j$ with at least four crossings.}%
        \label{fig:get-torus}
\end{figure}

\begin{sublem}\label{sublem:at-most-3}
We may assume that no two arcs $\bar l_i$ and $\bar l_j$ can meet
more than three times.
\end{sublem}

{\sc Proof:} Suppose the two arcs $\bar l_i$ and $\bar l_j$ meet
at least four times. We now choose an orientation of the knot
$K$. This orientation induces an orientation on $\pi(K)$ and its
strings $\bar l_i$ and $\bar l_j$. Now there are two cases to
consider:

\noindent{\sc Case 1. \sl The orientations of $\bar l_i$ and
$\bar l_j$ are consistent around the boundary of
$\delta_i\cap\delta_j$.}\ --- Notice that $\bar l_i$ and $\bar
l_j$ create at least four consecutive crescents such that one end
of the first is at $\pi(\mathcal Q)$ and no ends of the last is
at $\pi(\mathcal Q)$. Straightening the two complementary subarcs
of the two subarcs around the four crescents between $\bar l_i$
and $\bar l_j$, we are able to obtain a five crossing knot $T$ as
depicted in Figure~\ref{fig:get-torus}.
Sublemma~\ref{sublem:crescent} guarantees that $T$ is alternating
and hence is a torus knot of type $(2,5)$ which has superbridge
index 4. By Lemma~\ref{lem:straighten}, we get a contradiction
$4\le s(T)\le s(K)=3$.

\noindent{\sc Case 2. \sl The orientations of  $\bar l_i$ and
$\bar l_j$ are inconsistent around the boundary of
$\delta_i\cap\delta_j$.}\ --- In this case, only three crossings
of $\bar l_i\cup\bar l_j$ are required to draw a contradiction.

For any crossing point $Z$ of $\bar l_i\cup\bar l_j$, let $Z_0$
be the middle point of the two points $Z_i=\pi^{-1}(Z)\cap l_i$
and $Z_j=\pi^{-1}(Z)\cap l_j$. Let $O$ be a point on $\mathcal Q$
located so as to separate the four points of $K\cap\mathcal Q$
two by two.
This case breaks into two subcases:

\smallskip\noindent{\sc Subcase 2.1. \sl
The closure of $l_i\cup l_j$ is not connected.}\
--- Let $A$ and $B$ denote the starting point and the end point of the
oriented arc $l_i$, respectively, and let $C$ and $D$ denote the
starting point and the end point of the oriented arc $l_j$,
respectively. Notice that $A,B,C,D$ are the four points of
$K\cap\mathcal Q$.

Suppose that $O$ separates $A$ and $D$. Then it also separates
$B$ and $C$. Let $X$ and $Y$ be the first and second crossing of
$\bar l_i\cup\bar l_j$ along $\bar l_i$, respectively. Then the
three parallel lines $\mathcal Q$, $\pi^{-1}(X)$ and
$\pi^{-1}(Y)$ cut $K$ into eight disjoint arcs. Consider the
plane $\mathcal E$ determined by the three points $O$, $X_0$ and
$Y_0$. Then each of the following six arcs, two from $K-(l_i\cup
l_j)$, the two between $\mathcal Q$ and $\pi^{-1}(X)$, and the
two between $\pi^{-1}(X)$ and $\pi^{-1}(Y)$, crosses $\mathcal E$
because the end points are separated by the plane. If $O$
separates $A$ and $B$, then it also separates $C$ and $D$. In
this case, each of the remaining two arcs between $\pi^{-1}(Y)$
and $\mathcal Q$ crosses $\mathcal E$. If $O$ does not separate
$A$ and $B$ then it does not separate $C$ and $D$ either. In this
case, the existence of the crossing next to $Y$ guarantees that
the union of the two remaining arcs between $\pi^{-1}(Y)$ and
$\mathcal Q$ crosses $\mathcal E$ at least twice. Consequently,
the knot $K$ crosses $\mathcal E$ at least eight times, resulting
a contradiction $s(K)\ge4$.

Suppose that $O$ does not separate $A$ and $D$. Then it does not
separate $B$ and $C$ either. Then, according to
Sublemma~\ref{sublem:only-2-kinds}, one of the two arcs of
$K-(l_i\cup l_j)$ which corresponds to a simple loop in
Figure~\ref{fig:only-2-kinds}, together with the segment on
$\mathcal Q$ between its end points,  must bound an embedded disk
whose interior does not meet $K$. According to
\cite[Lemma~13]{kuperberg}, this kind of {\em topological
triviality\/} can be avoided at the beginning when we choose the
quadrisecant $\mathcal Q$. We now assume that our quadrisecant
$\mathcal Q$ is topologically nontrivial.

The quadrisecant $\mathcal Q$ of a nontrivial knot $K$ is defined
to be {\em topologically nontrivial\/} if, for any two points
$P_1$, $P_2$ of $K\cap\mathcal Q$ which are adjacent along
$\mathcal Q$,  any disk (possibly singular) bounded by the line
segment $P_1P_2$ and the arc of $K-\mathcal Q$ whose end points
are $P_1$ and $P_2$ meets $K$ in its interior.

\smallskip\noindent{\sc Subcase 2.2. \sl
The closure of $l_i\cup l_j$ is connected.}\
---  We may assume the starting point of $l_j$ is the end point of $l_i$.
Let $A$ and $B$ denote the starting point and the end point of the
oriented arc $l_i$, respectively, and $C$ the the end point of
the oriented arc $l_j$. Notice that $A,B,C$ are three  points of
$K\cap\mathcal Q$. Let $D$ be the remaining point of
$K\cap\mathcal Q$.

Suppose $O$ separates $B$ from the two points $A$ and $C$. Then
$B$ and $D$ are on the same half-line of $\mathcal Q-O$. Let $X$
and $Y$ be the first and second crossing point of $\bar
l_i\cup\bar l_j$ along $\bar l_i$. We choose two points
$X_+\in\pi^{-1}(X)$ and $Y_-\in\pi^{-1}(Y)$ so that $$
\overset\longrightarrow{OX_+}=
\overset\longrightarrow{OX_0}+\frac{\|\overset\longrightarrow{X_iX_j}\|}%
                                   {\|\overset\longrightarrow{OA}\|}
                              \overset\longrightarrow{OA}
\mbox{\quad and\quad } \overset\longrightarrow{OY_-}=
\overset\longrightarrow{OY_0}+\frac{\|\overset\longrightarrow{Y_iY_j}\|}%
                                   {\|\overset\longrightarrow{OB}\|}
                              \overset\longrightarrow{OB}.
$$ Let $\mathcal E$ be the plane determined by the three points
$O$, $X_+$ and $Y_-$. Then each of the eight arcs in
$K-(\pi^{-1}(X)\cup\pi^{-1}(Y)\cup\mathcal Q)$ crosses $\mathcal
E$ because the end points are separated by the plane. So we get a
contradiction $s(K)\ge4$.

Suppose $O$ separates $A$ from the two points $B$ and $C$. Then
$A$ and $D$ are on the same half-line of $\mathcal Q-O$. Let
$l_a$ be the component of $K-\mathcal Q$ whose end points are at
$A$ and $D$, and let $l_b$ be the one whose end points are at $C$
and $D$. By the assumption that $\mathcal Q$ is topologically
nontrivial, we know that $\bar l_b$ is the only arc corresponding
to a simple loop of Figure~\ref{fig:only-2-kinds}. Therefore
there is a crossing point $Y$ between $\bar l_a$ and $\bar
l_i\cup\bar l_j$. Consider the simple loop in $\tilde
l_i\cup\tilde l_j$ created by the last crossing point of $\bar
l_i\cup\bar l_j$ along $\bar l_i$. By
Sublemma~\ref{sublem:crescent}, this loop must have a crossing
with $\bar l_a$. Again, by the assumption that $\mathcal Q$ is
topologically nontrivial, $Y$ can be chosen so that the two
vectors $\overset\longrightarrow{Y_0Y_a}$ and
$\overset\longrightarrow{OA}$ are in opposite directions. Let $X$
be the first crossing point of $\bar l_i\cup\bar l_j$ along $\bar
l_i$ and let $\mathcal E$ be the plane determined by the three
points $O$, $X_0$ and $Y_0$. Again we consider the eight disjoint
arcs of $K-(\pi^{-1}(X)\cup\pi^{-1}(Y)\cup\mathcal Q)$. Each of
the following five, two from $l_a-\pi^{-1}(Y)$, the arc $l_b$,
and the two between $\mathcal Q$ and $\pi^{-1}(X)$, crosses the
plane $\mathcal E$. It remains to check how many times the
remaining three arcs cross $\mathcal E$. If $\bar l_a$ crosses
$\bar l_i$ at $Y$, the three arcs joins the points $X_i$, $Y_i$,
$B$, and $X_j$, successively. In this case, each of the three
arcs crosses $\mathcal E$. If $\bar l_a$ crosses $\bar l_j$ at
$Y$, the three arcs joins the points $X_i$, $B$, $Y_j$ and $X_j$,
successively. In this case, the arc joining $B$ and $Y_j$ crosses
$\mathcal E$ and the existence of a crossing point in $\bar
l_i\cup\bar l_j$ other than $X$ guarantees that the union of the
two remaining arcs crosses the plane at least twice.
Consequently, $K$ crosses $\mathcal E$ at least eight times,
resulting a contradiction $s(K)\ge4$.

The case when $O$ separates $C$ from the two points $A$ and $C$
can be handled similarly.~\qed

\begin{figure}[h]
\newcounter{PatternNum}
\setcounter{PatternNum}{1}
\def\patnum{\put(-28,15){\scriptsize(\arabic{PatternNum})}
            \addtocounter{PatternNum}{1}}
\begin{picture}(50,50)(-25,-22)
\patnum \put(0,0){\circle{40}} \qbezier(20,0)(0,0)(14.14,14.14)
\qbezier(0,20)(0,0)(-14.14,14.14)
\qbezier(-20,0)(0,0)(-14.14,-14.14)
\qbezier(0,-20)(0,0)(14.14,-14.14)
\end{picture}
\begin{picture}(50,50)(-25,-22)
\patnum \put(0,0){\circle{40}} \qbezier(20,0)(0,0)(14.14,14.14)
\qbezier(0,20)(0,0)(-14.14,14.14)
\qbezier(-20,0)(0,0)(14.14,-14.14)
\qbezier(-14.14,-14.14)(0,0)(0,-20)
\end{picture}
\begin{picture}(50,50)(-25,-22)
\patnum \put(0,0){\circle{40}} \qbezier(20,0)(0,0)(14.14,14.14)
\qbezier(0,20)(0,0)(-14.14,14.14) \put(-20,0){\line(1,-1){20}}
\put(-14.14,-14.14){\line(1,0){28.285}}
\end{picture}
\begin{picture}(50,50)(-25,-22)
\patnum \put(0,0){\circle{40}} \qbezier(20,0)(0,0)(14.14,14.14)
\put(0,20){\line(0,-1){40}} \qbezier(-14.14,14.14)(0,0)(-20,0)
\put(-14.14,-14.14){\line(1,0){28.285}}
\end{picture}
\begin{picture}(50,50)(-25,-22)
\patnum \put(0,0){\circle{40}} \qbezier(20,0)(0,0)(14.14,14.14)
\qbezier(0,20)(0,0)(14.14,-14.14)
\qbezier(-14.14,14.14)(0,0)(0,-20)
\qbezier(-20,0)(0,0)(-14.14,-14.14)
\end{picture}
\begin{picture}(50,50)(-25,-22)
\patnum \put(0,0){\circle{40}} \qbezier(20,0)(0,0)(14.14,14.14)
\put(0,20){\line(0,-1){40}}
\put(-14.14,14.14){\line(1,-1){28.285}}
\qbezier(-20,0)(0,0)(-14.14,-14.14)
\end{picture}

\begin{picture}(50,50)(-25,-22)
\patnum \put(0,0){\circle{40}} \qbezier(20,0)(0,0)(14.14,14.14)
\put(0,20){\line(-1,-1){20}} \qbezier(-14.14,14.14)(0,0)(0,-20)
\put(-14.14,-14.14){\line(1,0){28.285}}
\end{picture}
\begin{picture}(50,50)(-25,-22)
\patnum \put(0,0){\circle{40}} \qbezier(20,0)(0,0)(14.14,14.14)
\put(0,20){\line(0,-1){40}}
\put(-14.14,14.14){\line(0,-1){28.285}}
\qbezier(-20,0)(0,0)(14.14,-14.14)
\end{picture}
\begin{picture}(50,50)(-25,-22)
\patnum \put(0,0){\circle{40}} \qbezier(20,0)(0,0)(14.14,14.14)
\qbezier(0,20)(0,0)(-14.14,-14.14)
\put(-14.14,14.14){\line(1,-1){28.285}}
\put(-20,0){\line(1,-1){20}}
\end{picture}
\begin{picture}(50,50)(-25,-22)
\patnum \put(0,0){\circle{40}} \qbezier(20,0)(0,0)(14.14,14.14)
\qbezier(0,20)(0,0)(-14.14,-14.14)
\qbezier(-14.14,14.14)(0,0)(0,-20)
\qbezier(-20,0)(0,0)(14.14,-14.14)
\end{picture}
\begin{picture}(50,50)(-25,-22)
\patnum \put(0,0){\circle{40}} \qbezier(20,0)(0,0)(14.14,14.14)
\qbezier(0,20)(0,0)(14.14,-14.14)
\put(-14.14,14.14){\line(0,-1){28.285}}
\put(-20,0){\line(1,-1){20}}
\end{picture}
\begin{picture}(50,50)(-25,-22)
\patnum \put(0,0){\circle{40}} \put(20,0){\line(-1,1){20}}
\put(14.14,14.14){\line(-1,0){28.285}}
\put(-20,0){\line(1,-1){20}}
\put(-14.14,-14.14){\line(1,0){28.285}}
\end{picture}

\begin{picture}(50,50)(-25,-22)
\patnum \put(0,0){\circle{40}} \put(20,0){\line(-1,1){20}}
\qbezier(14.14,14.14)(0,0)(-20,0)
\qbezier(-14.14,14.14)(0,0)(0,-20)
\put(-14.14,-14.14){\line(1,0){28.285}}
\end{picture}
\begin{picture}(50,50)(-25,-22)
\patnum \put(0,0){\circle{40}} \put(20,0){\line(-1,1){20}}
\put(14.14,14.14){\line(-1,-1){28.285}}
\put(-14.14,14.14){\line(1,-1){28.285}}
\put(-20,0){\line(1,-1){20}}
\end{picture}
\begin{picture}(50,50)(-25,-22)
\patnum \put(0,0){\circle{40}} \put(20,0){\line(-1,1){20}}
\put(14.14,14.14){\line(-1,-1){28.285}}
\qbezier(-14.14,14.14)(0,0)(0,-20)
\qbezier(-20,0)(0,0)(14.14,-14.14)
\end{picture}
\begin{picture}(50,50)(-25,-22)
\patnum \put(0,0){\circle{40}} \qbezier(20,0)(0,0)(-14.14,14.14)
\qbezier(14.14,14.14)(0,0)(0,-20)
\qbezier(0,20)(0,0)(-14.14,-14.14)
\qbezier(-20,0)(0,0)(14.14,-14.14)
\end{picture}
\begin{picture}(50,50)(-25,-22)
\patnum \put(0,0){\circle{40}} \qbezier(20,0)(0,0)(-14.14,14.14)
\put(14.14,14.14){\line(-1,-1){28.285}}
\put(0,20){\line(0,-1){40}} \qbezier(-20,0)(0,0)(14.14,-14.14)
\end{picture}
\begin{picture}(50,50)(-25,-22)
\patnum \put(0,0){\circle{40}} \put(20,0){\line(-1,0){40}}
\put(14.14,14.14){\line(-1,-1){28.285}}
\put(0,20){\line(0,-1){40}}
\put(-14.14,14.14){\line(1,-1){28.285}}
\end{picture}
\caption{Patterns near $\pi(\mathcal Q)$}\label{fig:near-Q}
\end{figure}
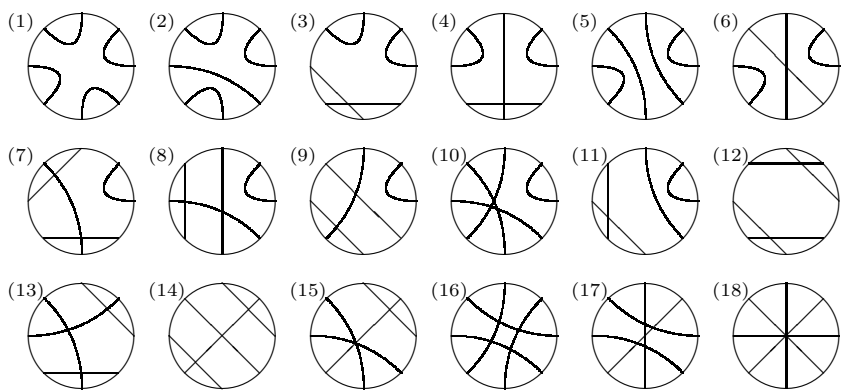

\begin{sublem}\label{sublem:only-finite}
There are finitely many possible diagrams for $K$ obtained from
$\pi(K)$ by perturbing near the quadrisecant $\mathcal Q$.
\end{sublem}

{\sc Proof:} Sublemma~\ref{sublem:only-2-kinds} and
Sublemma~\ref{sublem:at-most-3} leave only finitely many possible
projections of $K$ on $\mathcal Q^\perp$ outside a small
neighborhood of $\pi(\mathcal Q)$ up to planar isotopies. On the
other hand, the projection of a cylindrical neighborhood of a
quadrisecant of an arbitrary knot has eighteen possible patterns
as shown in Figure~\ref{fig:near-Q} where all four arcs are
smooth, up to small perturbations and planar isotopies. For each
pair of the projections outside and inside a neighborhood of
$\pi(\mathcal Q)$, there are only finitely many ways to combine
them to obtain a projection whose singular points are only
transversal double points. For each double point, there are only
two choices for crossings. Consequently there are only finitely
many possible diagrams of $K$ on $\mathcal Q^\perp$.~\qed

\medskip Since there are only finitely many possible diagrams, there are only finitely many possible knot types.

\section{An example}
\begin{figure}[t]
\includegraphics[210pt,307pt][388pt,451pt]{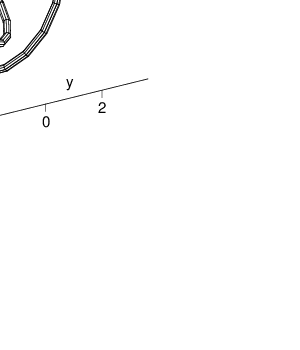}
\includegraphics[228pt,330pt][385pt,470pt]{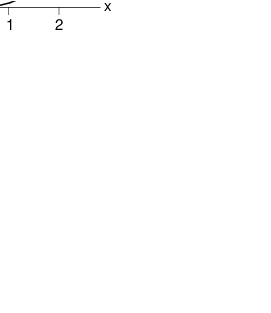}

(a)\hskip2.2in(b) \caption{A figure eight knot with harmonic
degree 3}\label{fig:fig8-aaron}
\end{figure}

\begin{figure}[t]
\includegraphics[203pt,314pt][393pt,446pt]{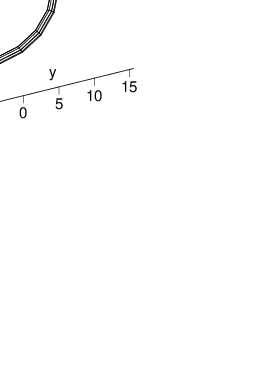}
\includegraphics[229pt,325pt][385pt,480pt]{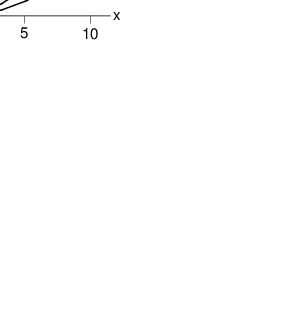}

(a)\hskip2.2in(b)
\caption{A figure eight knot having $z$-axis as its quadrisecant}%
        \label{fig:fig8-taek}
\end{figure}

The knot shown in Figure~\ref{fig:fig8-aaron}(a) is a figure
eight knot parametrized\footnote{This parametrization was
obtained by modifying Trautwein's parametrization
in~\cite{trautwein}:\\
\def\mytab{\hbox to2cm{\hfill}}
\mytab$x(t)=32\cos t-51\sin t-104\cos2t-34\sin2t+104\cos3t-
91\sin3t$\\ \mytab$y(t)=94\cos t+41\sin t+113\cos2t         -
68\cos3t-124\sin3t$\\ \mytab$z(t)=16\cos t+73\sin
t-211\cos2t-39\sin2t- 99\cos3t- 21\sin3t$ } by
\begin{eqnarray*}
x(t)&=&307\cos^3t+5346\sin t\cos^2t-2663\cos^2t\\
     &&\hskip0.5in{}-26\sin t\cos t-1142\cos t-1378\sin t+1280\\
y(t)&=&6337\cos^3t+191\sin t\cos^2t+691\cos^2t\\
     &&\hskip0.5in{}+103\sin t\cos t-5021\cos t-1019\sin t+677\\
z(t)&=&373\cos^3t-3157\sin t\cos^2t-4436\cos^2t\\
     &&\hskip0.5in{}-1029\sin t\cos t+50\cos t+910\sin t+2222
\end{eqnarray*}
for $0\le t\le2\pi$. Since its harmonic degree~\cite{trautwein}
is 3, its superbridge number is 3. Figure~\ref{fig:fig8-aaron}(b)
shows its projection into the $xy$-plane. Up to scaling and
reparametrization, this knot can be perturbed to have a
polynomial parametrization
\begin{eqnarray*}
x(t)&=&(2t-1)(4t-1)(10t-1)(25t-16)(25t-21)(50t-9)\times\\
&&(386t^6-708t^5-201t^4+945t^3-383t^2
-\frac{42224361}{1146679}t-\frac{2701080}{1146679})\\
y(t)&=&-70(2t-1)^2(4t-1)(10t-1)(25t-21)^2\times\\
&&(229t^6-776t^5+806t^4-197t^3-56 t^2
-\frac{1667040}{277477}t-\frac{104544}{277477})\\
z(t)&=&(20t-3)(25t-9)(25t-16)(25t-23)(1233t^8-5985t^7+11394t^6\\
&&\hskip0.5in{}-10375t^5+4167t^4-243t^3-179t^2
-\frac{2145804}{166595}t-\frac{712368}{832975})
\end{eqnarray*}
for $0\le t\le1$. Figure~\ref{fig:fig8-taek} illustrates this
knot. It is clear that the $z$-axis is a quadrisecant of this
knot.\footnote{ This knot may have superbridge number bigger than
3. However, it can be reduced to 3 again, by applying
Lemma~\ref{lem:straighten} away from the $z$-axis.} It is a
combination of Figure~\ref{fig:only-2-kinds}(b) and
Figure~\ref{fig:near-Q}(17). Figure~\ref{fig:fig8-aaron} and
Figure~\ref{fig:fig8-taek} are of scale $1:1000$.

\section{Prime knots up to 9 crossings}
All 3-superbridge knots are among the forty seven knots which are
2-bridge knots up to 9 crossings except the torus knots of types
$(2,5)$, $(2,7)$ and $(2,9)$. They are marked with $\star$ or
$\times$ in Table~\ref{tab:1}. The symbols in the first column
are as in~\cite{burde-zieschang,rolfsen}. The number 47 is a very
rough upper bound for the number of 3-superbridge knots. To show
that it is an upper bound, we only need to show that
3-superbridge knots cannot have minimal crossing number bigger
than 9.

\begin{figure}[t]
\small
\begin{picture}(135, 97.5)(-67.5,-18.75)
\put(-5,-26){(a)} \put(-45,67.5){\line(1,0){30}}
\put(-15,67.5){\line(0,-1){22.5}} \put(-15,45){\line(-1,0){30}}
\put(-45,45){\line(0,1){22.5}} \put(15,67.5){\line(1,0){30}}
\put(45,67.5){\line(0,-1){22.5}} \put(45,45){\line(-1,0){30}}
\put(15,45){\line(0,1){22.5}} \put(-15,60){\line(1,0){30}}
\qbezier[27](45,60)(67.5,60)(67.5,37.5)
\qbezier[60](67.5,37.5)(67.5,7.5)(0,0)
\qbezier(0,0)(-60,22.5)(-60,37.5)
\qbezier(-60,37.5)(-60,52.5)(-45,52.5)
\qbezier[30](-15,52.5)(-7.5,52.5)(0,0) \thicklines
\qbezier(-15,52.5)(0,52.5)(0,63.75)
            \qbezier(0,63.75)(0,75)(15,75)
            \put(15,75){\line(1,0){30}}
            \qbezier(45,75)(52.5,75)(52.5,67.5)
            \qbezier(52.5,67.5)(52.5,60)(45,60)\thinlines
\qbezier(0,0)(7.5,-15)(0,-15) \qbezier(0,-15)(-7.5,-15)(0,0)
\qbezier(0,0)(7.5,52.5)(15,52.5)
\qbezier(45,52.5)(60,52.5)(60,37.5)
\qbezier(60,37.5)(60,22.5)(0,0)
\qbezier(0,0)(-67.5,7.5)(-67.5,37.5)
\qbezier(-67.5,37.5)(-67.5,60)(-45,60)
\end{picture}
\quad
\begin{picture}(135, 97.5)(-67.5,-18.75)
\put(-6,-26){(b)} \put(-45,67.5){\line(1,0){30}}
\put(-15,67.5){\line(0,-1){22.5}} \put(-15,45){\line(-1,0){30}}
\put(-45,45){\line(0,1){22.5}} \put(15,67.5){\line(1,0){30}}
\put(45,67.5){\line(0,-1){22.5}} \put(45,45){\line(-1,0){30}}
\put(15,45){\line(0,1){22.5}} \put(-15,60){\line(1,0){30}}
\qbezier(45,60)(67.5,60)(67.5,37.5)
\qbezier(67.5,37.5)(67.5,7.5)(0,0)
\qbezier(0,0)(-60,22.5)(-60,37.5)
\qbezier(-60,37.5)(-60,52.5)(-45,52.5)
\qbezier[30](-15,52.5)(-7.5,52.5)(0,0) \thicklines
\qbezier(-15,52.5)(0,52.5)(0,63.75) \qbezier(0,63.75)(0,75)(15,75)
\qbezier(15,75)(75,75)(75,37.5) \qbezier(75,37.5)(75,-15)(0,-15)
\thinlines \qbezier[11](0,0)(7.5,-15)(0,-15)
\qbezier(0,-15)(-7.5,-15)(0,0) \qbezier(0,0)(7.5,52.5)(15,52.5)
\qbezier(45,52.5)(60,52.5)(60,37.5)
\qbezier(60,37.5)(60,22.5)(0,0)
\qbezier(0,0)(-67.5,7.5)(-67.5,37.5)
\qbezier(-67.5,37.5)(-67.5,60)(-45,60)
\end{picture}
\caption{Moves reducing the number of
crossings}\label{fig:at-most-9}
\end{figure}
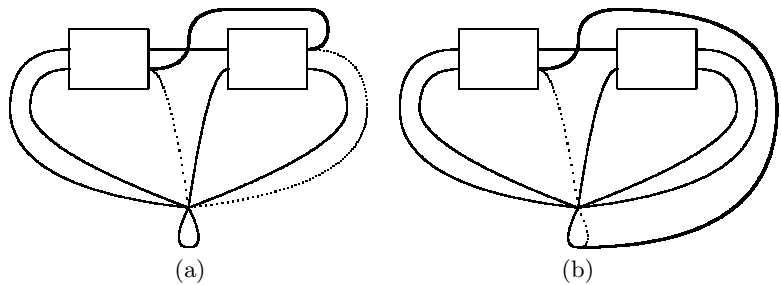

By Sublemma~\ref{sublem:at-most-3}, we know that there are at most
three double points in Figure~\ref{fig:only-2-kinds}(a) and at
most six in Figure~\ref{fig:only-2-kinds}(b). After a little
perturbation if necessary, each pattern in
Figure~\ref{fig:near-Q} have at most six crossings. Therefore
knot diagrams obtained from Figure~\ref{fig:only-2-kinds}(a)
cannot have more than nine crossings. Since each pattern in
Figure~\ref{fig:near-Q}(1)--(14) has at most three crossings,
knot diagrams obtained by any combination of
Figure~\ref{fig:only-2-kinds}(b) and one of
Figure~\ref{fig:near-Q}(1)--(14) cannot have more than nine
crossings. It remains to handle the combinations of
Figure~\ref{fig:only-2-kinds}(b) and one of
Figure~\ref{fig:near-Q}(15)--(18). Since the quadrisecant
$\mathcal Q$ meets the knot $K$ at four distinct points, the four
arcs in any of Figure~\ref{fig:near-Q}(15)--(18) are in distinct
vertical levels. Therefore the crossings obtained from
Figure~\ref{fig:near-Q}(15)--(18) are not alternating, and hence
any combination of Figure~\ref{fig:only-2-kinds}(b) and one of
Figure~\ref{fig:near-Q}(15)--(16) gives a non-alternating diagram
of at most ten crossings. Because 3-superbridge knots are
alternating knots, those obtained from such combinations must
have minimal crossing number at most 9. Now it remains to
consider the combinations of Figure~\ref{fig:only-2-kinds}(b) and
one of Figure~\ref{fig:near-Q}(17)--(18). For any such
combinations, we are able to move the uppermost arc or the
lowermost arc at the quadruple point to reduce the number of
crossings by one or two as needed to reduce the number down to at
most 10. Two examples of such moves are shown in
Figure~\ref{fig:at-most-9}. The resulting diagrams are still
nonalternating, hence their minimal crossing numbers are at most
9.

\def\tableone{
\begin{tabular}{|c|c|c|}
\hline $K$&$s[K]$&\\ \hline $3_1$& 3&$\circ$$\star$\\ $4_1$&
3&$\star$\\ $5_1$& 4&$\circ$\\ $5_2$& 3--4&$\star$\\ $6_1$&
3--4&$\star$\\ $6_2$& 3--4&$\star$\\ $6_3$& 3--4&$\star$\\ $7_1$&
4&$\circ$\\ $7_2$& 3--4&$\star$\\ $7_3$& 3--4&$\star$\\ $7_4$&
3--4&$\star$\\ $7_5$& 4&$\times$\\ $7_6$& 4&$\times$\\ $7_7$&
4&$\times$\\ $8_1$& 4--5&$\times$\\ $8_2$& 4--5&$\times$\\ $8_3$&
4--6&$\times$\\ \hline
\end{tabular}}
\def\tabletwo{
\begin{tabular}{|c|c|c|}
\hline $K$&$s[K]$&\\ \hline $8_4$& 3--5&$\star$\\ $8_5$&
4--6&$\diamond$\\ $8_6$& 4--6&$\times$\\ $8_7$& 3--6&$\star$\\
$8_8$& 4--5&$\times$\\ $8_9$& 3--6&$\star$\\ $8_{10}$&
4--6&$\diamond$\\ $8_{11}$& 4--5&$\times$\\ $8_{12}$&
4--6&$\times$\\ $8_{13}$& 4--5&$\times$\\ $8_{14}$&
4--5&$\times$\\ $8_{15}$& 4--6&$\diamond$\\ $8_{16}$&
4&$\diamond$\\ $8_{17}$& 4&$\diamond$\\ $8_{18}$& 4&$\diamond$\\
$8_{19}$& 4&$\circ$$\diamond$\\ $8_{20}$& 4&$\diamond$\\ \hline
\end{tabular}}
\def\tablethree{
\begin{tabular}{|c|c|c|}
\hline $K$&$s[K]$&\\ \hline $8_{21}$& 4&$\diamond$\\ $9_1$&
4&$\circ$\\ $9_2$& 4--7&$\times$\\ $9_3$& 4--6&$\times$\\ $9_4$&
4--7&$\times$\\ $9_5$& 4--6&$\times$\\ $9_6$& 4--6&$\times$\\
$9_7$& 4--6&$\times$\\ $9_8$& 4--6&$\times$\\ $9_9$&
4--6&$\times$\\ $9_{10}$& 4--6&$\times$\\ $9_{11}$&
4--6&$\times$\\ $9_{12}$& 4--6&$\times$\\ $9_{13}$&
4--6&$\times$\\ $9_{14}$& 4--7&$\times$\\ $9_{15}$&
4--5&$\times$\\ $9_{16}$& 4--7&$\diamond$\\ \hline
\end{tabular}}
\def\tablefour{
\begin{tabular}{|c|c|c|}
\hline $K$&$s[K]$&\\ \hline $9_{17}$& 4--7&$\times$\\ $9_{18}$&
4--6&$\times$\\ $9_{19}$& 4--6&$\times$\\ $9_{20}$&
4--6&$\times$\\ $9_{21}$& 4--7&$\times$\\ $9_{22}$&
4--7&$\diamond$\\ $9_{23}$& 4--7&$\times$\\ $9_{24}$&
4--6&$\diamond$\\ $9_{25}$& 4--7&$\diamond$\\ $9_{26}$&
4--6&$\times$\\ $9_{27}$& 4--6&$\times$\\ $9_{28}$&
4--6&$\diamond$\\ $9_{29}$& 4--7&$\diamond$\\ $9_{30}$&
4--6&$\diamond$\\ $9_{31}$& 4--6&$\times$\\ $9_{32}$&
4--6&$\diamond$\\ $9_{33}$& 4--6&$\diamond$\\ \hline
\end{tabular}}
\def\tablefive{
\begin{tabular}{|c|c|c|}
\hline $K$&$s[K]$&\\ \hline $9_{34}$& 4--6&$\diamond$\\ $9_{35}$&
4--6&$\diamond$\\ $9_{36}$& 4--7&$\diamond$\\ $9_{37}$&
4--7&$\diamond$\\ $9_{38}$& 4--7&$\diamond$\\ $9_{39}$&
4--6&$\diamond$\\ $9_{40}$& 4&$\diamond$\\ $9_{41}$&
4&$\diamond$\\ $9_{42}$& 4&$\diamond$\\ $9_{43}$&
4--5&$\diamond$\\ $9_{44}$& 4--5&$\diamond$\\ $9_{45}$&
4--5&$\diamond$\\ $9_{46}$& 4&$\diamond$\\ $9_{47}$&
4--6&$\diamond$\\ $9_{48}$& 4--6&$\diamond$\\ $9_{49}$&
4--5&$\diamond$\\
        &     & \\
\hline
\end{tabular}}
\begin{table}[b]
\small \tableone\tabletwo\tablethree\tablefour\tablefive
\bigskip
\caption{Superbridge index of prime knots up to 9 crossings}%
       \label{tab:1}
\end{table}

Among all possible combinations of one of
Figure~\ref{fig:only-2-kinds} and one of Figure~\ref{fig:near-Q},
we do not get all the forty seven knots mentioned above. There
are 35 knots which either do not appear during the construction or
are excluded by using the methods used in the proof of
Theorem~\ref{thm:3-sup}.  They are marked with $\times$ in the
table. The remaining 12 knots are marked with $\star$. This list
of 12 knots contains all the 3-superbridge knots. Details behind
the selection of 12 knots and the rejection of 35 knots will be
handled in~\cite{jeon-jin}.

In the table, torus knots are marked with $\circ$, for which the
superbridge index is determined by Theorem~\ref{thm:kuiper-torus}.

If a knot is presented as a polygon in space, one half of the
number of edges is an upper bound of the superbridge
index~\cite{jin-poly}. The number or the upper limit of the range
of numbers in the second column of the table is the largest
integer not exceeding one half of the minimal edge number or the
best-known minimal edge
number~\cite{calvo,calvo-millett,jin-poly,meissen}. For the five
knots, $3_1$, $4_1$, $5_1$, $5_2$ and $8_{19}$, this number is
equal to the harmonic degree~\cite{trautwein}. It is known that
2-bridge knots cannot have superbridge index bigger than
seven~\cite{furstenberg-et-al}.

On the other hand, for those marked with $\star$ or $\diamond$,
the number or the lower limit of the range of numbers in the
second column is one bigger than the bridge index. For those
marked with $\times$ or only with $\circ$, the number or the
lower limit of the range of numbers in the second column is two
bigger than the bridge index.

Among the 18 knots whose superbridge index is determined in the
table, only $7_5$, $7_6$ and $7_7$ were newly found by this work.

\end{document}